\documentclass[a4paper,11pt]{article}
\usepackage{graphicx}
\usepackage{amsmath}
\usepackage{amssymb}
\usepackage{amstext}
\usepackage{color}
\usepackage{colortbl}
\usepackage{xcolor}

\newcommand{\proof}{\vskip 5truemm \noindent{\textsc Proof.}~}

\newtheorem{thm}{Theorem}[section]

\newtheorem{lem}{Lemma}[section]

\numberwithin{equation}{section}
\begin{document}

%\title{\bf Existence of right continuous with bounded variation solution of a perturbed maximal monotone
%evolution}
\title{\bf A measure differential inclusion involving time-dependent maximal monotone operators}

\author{Dalila Azzam-Laouir
 \footnote{LAOTI, FSEI, Universit\'e Mohamed Seddik Benyahia de Jijel, Alg\'erie. E-mail: laouir.dalila@gmail.com,
 dalilalaouir@univ-jijel.dz}}
%\hskip 4pt Charles Castaing\footnote{  C. Castaing IMAG, Univ
%Montpellier, CNRS, Montpellier II, 34095, Case courrier 051,
%Montpellier Cedex 5, France. E-Mail: charles.castaing@gmail.com}
%\hskip 4pt M. D. P. Monteiro Marques \footnote{CMAF and Faculdade de
%Ciencias de Lisboa, Av. Prof. Gama Pinto 2, P. 1600 Lisboa,
%Portugal. E-Mail: mmarques@lmc.fc.ul.pt}}

\maketitle
\begin{abstract} We establish in this paper, an existence and
uniqueness result of right continuous with bounded variation
solution for a perturbed differential inclusion governed by
time-dependent maximal monotone operators.

\end{abstract}
\vskip4mm \noindent{\footnotesize\textbf{Keywords:} Bounded
variation, differential measure, Lipschitz mapping,  maximal
monotone operator, pseudo-distance, right continuous.}

\vskip4mm \noindent{\footnotesize\textbf{AMS Subject
Classifications: 2010}: 34H05, 34K35, 60H10 28A25, 28C20}

\vskip6mm

%------------------------------------------------------------------------------
\section{Introduction}
Let $I=[0,T]$ ($T>0$). We consider in this paper, in a separable
Hilbert space $\mathcal{H}$, the following perturbed evolution
differential inclusion
\begin{equation}\label{1.1} -Du(t) \in A(t) u(t) +f(t,
u(t))\;\;a.e.,\;\;\;u(0)=u_0,\end{equation} where for each $t\in I$,
$A(t)$ is a maximal monotone operator of $\mathcal{H}$,  and where
the set-valued map $t \mapsto A(t) $ is right continuous with
bounded variation (BVRC),
  in the sense  that
 there exists a function  $\rho:  I \rightarrow [0, \infty[$,
which is right continuous on $[0, T[$  and nondecreasing  with
$\rho(0) = 0$ and  $\rho(T) < \infty$ such that
$$ dis(A(t), A(s)) \leq d\rho(]s, t]) = \rho(t)-\rho(s) \;\;\; 0\leq s \leq t \leq T,$$
here  $dis(\cdot, \cdot)$ is the pseudo-distance between maximal
monotone operators introduced by Vladimirov \cite{Vla}; see relation
\eqref{2.1}, and $f:I\times \mathcal{H}\longrightarrow \mathcal{H}$
is  measurable in $t$ and satisfies a Lipschitz condition w.r.t the
second variable. For the study of problem \eqref{1.1} without
perturbation, we refer to Theorem 1 in \cite{KM}.

Our result extends to the BVRC time-depend maximal monotone
operators, Theorem 3.1 in \cite{ACM}, dealing with BVC (continuous
with bounded variation) time-dependent maximal monotone operators,
and Theorem 4.1 in \cite{AHT} dealing with BVRC sweeping processes.

It is worth to mention that in the recent work by Azzam et al
\cite{ACM1}, the authors gave a theorem on the existence and
uniqueness of a BVRC solution to problem \eqref{1.1}, by assuming
that the domain of the operator $A(t)$ is ball-compact. In the
present paper, we will establish the same result without any
compactness assumption. Actually, we will prove the Cauchy criteria
to conclude the convergence of our approximate sequence to the
desired solution, by following  ideas in the proofs of Theorem 3.1
in \cite{ACM} and Theorem 4.1 in \cite{AHT}.

Let $\lambda$ be the Lebesgue measure on $I$ and $d\rho$  the
Stieljes measure associated with $\rho$.  We set $\nu := \lambda
+d\rho$ and $\frac{ d\lambda}{d\nu}$ the density of $ \lambda$ w.r.t
 $\nu$.  By a solution of \eqref{1.1}, we mean that
there exists a BVRC  mapping $u: I\longrightarrow \mathcal{H}$
satisfying
$$
(P_f)
\begin{cases}
u(0)=u_0\in D(A(0));\\
 u(t)\in D(A(t))\;\;\;\forall t\in I;\\
   -\displaystyle\frac{du}{d\nu }(t)\in A(t) u(t)+ f(t, u(t))   \frac{d\lambda} {d\nu}(t)  \;\;\;d\nu-a.e.\,t\in
   I,
\end{cases}
$$
where $\frac{du}{d\nu }$ denotes the density of $u$ w.r.t $\nu$.

There is an intensive study on differential inclusions governed by
time-dependent and time and state dependent maximal monotone
operators, where the variation of the operator is absolutely
continuous or continuous with bounded variation; see for instance
\cite{ABCM1, ABCM2, AB, ACM, BAC, Ken, KM, Le, P, SAM, TOL1, Vla}.
However, there are few works concerning  $(P_f)$ when the variation
of the operator $A(t)$ is right continuous with bounded variation
(BVRC). Actually, to the best of our knowledge, there is an
existence result to $(P_0)$ $(f=0)$ in \cite{KM} and for the
perturbed problem we refer to the recent reference \cite{ACM1}.

Differential inclusions governed by time-dependent maximal monotone
operators constitute generalization of the so called "sweeping
process", that is differential inclusions governed by the normal
cone to closed and convex moving sets, since this normal cone is a
maximal monotone operator. This process was introduced and widely
studied by J.J. Moreau; we cite for instance \cite{M, M1}. We can
also refer to the book of Monteiro-Marques \cite{Mon}. Following the
results of Moreau, several extensions in diverse directions, in
particular  perturbed sweeping processes, have been studied in the
literature, where the moving sets in the normal cone have absolutely
continuous variation or continuous and bounded in variation (BVC).
The BVRC case was also deeply discussed in many interesting papers;
we cite for instance \cite{AHT, ANT, BCSS, CM, ET1, M2, NNT, Thi3}.

The paper is organized as follows. In section 2, we give notations
and recall the preliminary results that we need in the sequel. In
section 3, we prove the existence and uniqueness  of  a right
continuous with bounded variation solution to the evolution problem
$(P_f)$ when $t \mapsto A(t) $ is BVRC and the perturbation $f:
I\times \mathcal{H}\longrightarrow \mathcal{H}$ is measurable on $I$
and satisfies a Lipschitz condition w.r.t the second variable.

\section{Notations and  Preliminaries} From now on, $I:=[0,T]$ $(T>0)$ is an interval of
$\mathbb{R}$ and $(\mathcal{H}, \langle \cdot ,\cdot \rangle)$ is a
separable Hilbert space. Its norm will  be denoted by $\|\cdot\|$,
and its unit closed ball (resp. closed ball of center $0$ and radius
$r>0$) will be denoted by $\overline{B}_{\mathcal{H}}$ (resp.
$r\overline{B}_{\mathcal{H}}$). We will denote by $\mathcal{B}(I)$
(resp. $\mathcal{B}(\mathcal{H})$) the Borel tribe on $I$ (resp. on
$\mathcal{H})$. The identity mapping of $\mathcal{H}$ will be
denoted by $Id_{\mathcal{H}}$. For a subset $K$ of $\mathcal{H}$,
$\overline{co}(K)$ will be the closed convex hull of $K$, which is
characterized by:
\begin{equation}\label{co}
\overline{co}(K)=\big\{y\in \mathcal{H}:\;\langle y, x\rangle\leq
\sup_{z\in K}\langle z, x\rangle\;\forall x\in K\big\}.
\end{equation}

If $\mu$ is a positive measure on $I$, we will denote by $L^p(I,
\mathcal{H}; \mu)$ $p \in [1, +\infty[$ (resp. $p=+\infty)$, the
Banach space of classes of $p$-$\mu$-integrable (resp.
$\mu$-essentially bounded) functions equipped with its classical
norm $\|\cdot\|_p$ (resp. $\|\cdot\|_{\infty}$).
%and we will denote
%by $\mathcal{C}(I,\mathcal{H})$ the Banach space of all continuous
%mappings $u : I\longrightarrow \mathcal{H}$, endowed with the sup
%norm.

We introduce in the following,  definitions and some properties of
functions with bounded variation and general vector measures. We
refer to \cite{Mon, M2, M3, MV} for more details.

Let $u:I\longrightarrow \mathcal{H}$. The variation of $u$
 in $I$ is the nonnegative extended real number
 \begin{equation}\label{1}
 var(u ; I):=\sup\sum_{k=1}^{n}\|u(t_k)-u(t_{k-1})\|,
 \end{equation}
 where the supremum is taken in $[0,+\infty]$ w.r.t all the finite
 sequences $t_0<t_1<\cdots<t_n$ of points of $I$ ($n$ is arbitrary).
 The function $u$ is said with bounded variation (BV) if and only if
 $var(u ; I)<+\infty$, and we write $u\in BV(I, \mathcal{H})$. In this
 case we have
 \begin{equation*}
 \lim_{s\uparrow t}\|u(t)-u(s)\|=\lim_{s\uparrow t} var(u; [s,t]).
 \end{equation*}

 Now, let $\mu$ be a positive Radon measure on $I$ and $\hat{\mu}$
 be an $\mathcal{H}$-valued measure on
 $I$ admitting a density $\frac{d \hat{\mu}}{d\mu}\in
 L^1_{loc}(I;\mathcal{H};\mu)$. Then for $d \mu$-almost every $t\in
 I$, we have
 \begin{equation}\label{2}
\frac{d \hat{\mu}}{d\mu}(t)=\lim_{\varepsilon\downarrow 0}\frac{d
\hat{\mu}([t, t+\varepsilon])}{d\mu([t,
t+\varepsilon])}=\lim_{\varepsilon\downarrow 0}\frac{d
\hat{\mu}([t-\varepsilon, t])}{d\mu([t-\varepsilon, t])}.
 \end{equation}
 The measure $\hat{\mu}$ is absolutely continuous w.r.t $\mu$ if and
 only if $\hat{\mu}=\frac{d \hat{\mu}}{d\mu}\mu$, i.e.,  $\frac{d
 \hat{\mu}}{d\mu}$ is a density of $\hat{\mu}$  w.r.t $\mu$. In this
 case a mapping $u:I\longrightarrow \mathcal{H}$ is
 $\hat{\mu}$-integrable if and only if the mapping $t\mapsto u(t)\frac{d
 \hat{\mu}}{d\mu}(t)$ is $\mu$-integrable, and we have
 \begin{equation}\label{7}
 \int_I u(t) d\hat{\mu}(t)=\int_I u(t)\frac{d
 \hat{\mu}}{d\mu}(t) d\mu(t).
 \end{equation}
If $u:I\longrightarrow \mathcal{H}$ is BVRC and $du$ is its
 differential measure, then we have
 \begin{equation}\label{3}
 u(t)=u(s)+\int_{]s,t]} du\;\;\;\;\;\forall s, t\in I\;(s\leq t).
 \end{equation}
 Conversely, if there exists $v\in L^1(I,\mathcal{H};\mu)$ such that
 $u(t)=u(0)+\int_{]0, t]} v \,d \mu$ for all $t\in I$, then $u$ is
 BVRC and $du=v \, d\mu$, that is $v$ is a density of the vector
 measure $du$ w.r.t the measure $\mu$. So that by \eqref{2}, we get
 \begin{equation}\label{4}
v(t)=\frac{d u}{d\mu}(t)=\lim_{\varepsilon\downarrow 0}\frac{ du
([t, t+\varepsilon])}{d\mu([t,
t+\varepsilon])}=\lim_{\varepsilon\downarrow 0}\frac{ du
([t-\varepsilon, t])}{d\mu([t-\varepsilon, t])}.
 \end{equation}
 If $\mu(\{t\})>0$, this last relation shows that
 \begin{equation}\label{5}
 v(t)=\frac{d u}{d\mu}(t)=\frac{d
 u(\{t\})}{d\mu(\{t\})}\;\;\;\textmd{and}\;\;\;\frac{d
 \lambda}{d\mu}(t)=0.
 \end{equation}

 The following lemmas will be needed in our proof.

 \begin{lem}\label{Lem2.1}(\cite{M3})
 Let $\mu$ be a positive Radon measure on $I$ and $u:I\longrightarrow
 \mathcal{H}$ a BVRC mapping such that its differential measure $du$ has a density
 $\frac{du}{d\mu}$ w.r.t $\mu$. Then the real function
 $\varphi:t\mapsto \varphi(t)=\|u(t)\|^2$ is BVRC and its
 differential measure $d\varphi$ satisfies
 \begin{equation}\label{6}
 d\varphi\leq 2\big\langle u, \frac{du}{d\mu}\big\rangle d\mu.
 \end{equation}
 \end{lem}

 \begin{lem}\label{lem2.6}(\cite{BCG}, Lemma 2.1)  Let $\mu$ be a positive Radon measure on $I$. Let $g \in L^1(I, \mathbb{R}; \mu)$ be a nonnegative function
 and $\beta  \geq 0$ be such that,
$\forall t \in I$, $0\leq  \mu ( \{t\} ) g(t)  \leq \beta < 1$. Let
$\varphi \in  L^\infty(I, \mathbb{R}; \mu)$ be a nonnegative
function satisfying
\begin{equation*}\varphi(t)  \leq \alpha +\int_{]0, t]}  g(s) \varphi(s)
\mu(ds)\;\;\;  \forall t \in I, \end{equation*} where $\alpha$ is a
nonnegative constant. Then
\begin{equation*}\varphi(t) \leq \alpha \exp \Big( \frac{1} {1-\beta} \int_{]0, t]} g(s) \mu(ds)\Big)\;\;\;    \forall t \in I. \end{equation*}
\end{lem}

%%%%%%%%%%%
%\begin{lem} \label{lem2.7} (\cite{ACM}, Lemma 2.7)  Let $\mu$ be a non-atomic positive Radon measure on the interval $I$.
%Let $c$, $p$ be nonnegative real functions
%such that $c \in L^1(I, \mathbb{R};\mu), p \in L^\infty (I,
%\mathbb{R};\mu)$, and let $\alpha\geq 0$. Assume that for
%all $ t \in I$,
%$\mu-a.e.\, \, t\in I$
%$$p(t) \leq \alpha +\int_{0}^{t} c(s) p(s) \mu(ds).$$
%Then, for $\mu-a.e.\, \, t\in I$
%$$p(t) \leq \alpha \exp\Big(\int _{0}^{t} c(s)
%\mu(ds)\Big).$$
%\end{lem}

%\begin{lem} \label{lem2.8} (Proposition 4.1 in \cite{TOL}) Let $m\in L^1(I, \mathbb{R};
%\lambda)$ be a nonnegative function,  $x:I\longrightarrow [0,
%+\infty[$ be a BVRC function and $a\geq 0$. If
%$$ \frac{1}{2} x^2(t)\leq \frac{1}{2} a^2+\int_{]0, t]} m(s)
%x(s)ds\;\;\;\forall t\in I,$$ then
%$$  x(t)\leq  a+\int_{]0, t]} m(s)
%ds\;\;\;\forall t\in I.$$
%\end{lem}

\begin{lem}\label{lem2.5}
Let $(\alpha_i)$, $(\beta_i)$, $(\gamma_i)$ and $(a_i)$ be sequences
of nonnegative real numbers such that $a_{i+1}\leq
\alpha_i+\beta_i\big(a_0+a_1+\cdots+a_{i-1}\big)+(1+\gamma_i)a_i$
for $i\in \mathbb{N}$. Then
$$ a_j\leq \bigg(a_0+\sum_{k=0}^{j-1}
\alpha_k\bigg)\exp\bigg(\sum_{k=0}^{j-1}
\big(k\beta_k+\gamma_k\big)\bigg)\;\;\textmd{for}\;j\in
\mathbb{N}^*.$$
\end{lem}

We finish this section by the definition and some properties of
maximal monotone operators. We refer the reader to \cite{Ba, brezis,
Vra} for these concepts.

 Let $A: \mathcal{H}
\rightrightarrows \mathcal{H}$ be a set-valued map. We denote by
$D(A)$, $R(A)$ and $Gr(A)$ its domain, range and graph. We say that
the operator $A$ is monotone, if $\langle y_1 -y_2, x_1 -
x_2\rangle\ge 0$ for all $(x_i, y_i)\in Gr(A)$ $(i=1, 2)$, and
  we say that $A$ is a maximal monotone operator of $\mathcal{H}$,
  if it is monotone and
its graph could not be contained strictly in the graph of any other
monotone operator.
%By Minty's Theorem, $A$ is maximal monotone iff
%$R(I_{\mathcal{H}}+A)=\mathcal{H}$.

If $A$ is a maximal monotone operator of $\mathcal{H}$, then  for
every $x \in D(A)$, $A(x)$ is nonempty closed and convex. We denote
the projection of the origin on the set $A(x)$ by $A^0(x)$.

 For $\eta > 0$, we denote by $J_{\eta}^A= (Id_{\mathcal{H}} +
\eta A)^{-1}$  the resolvent of  $A$. %and by $A_{\eta} = \frac{1}{
%\eta} (I_{\mathcal{H}}- J_{\eta}^A)$ its Yosida approximation.
It is
 well-known that
  this operator
is single-valued and defined on all of $\mathcal{H}$, furthermore
$J^A_{\eta}(x)\in D(A)$, for all $x\in \mathcal{H}$.

Let $A: D(A)\subset \mathcal{H}\rightrightarrows \mathcal{H}$ and
$B: D(B)\subset \mathcal{H}\rightrightarrows \mathcal{H}$ be two
maximal monotone operators. Then,  $dis(A, B)$ is the Vladimirov's
pseudo-distance between $A$ and $B$, defined as follows
\begin{equation}\label{2.1} dis(A,B)=\sup\bigg\{\frac{\langle y-y', x'-x\rangle}{1+\|
y\|+\| y'\|}:\;(x,y)\in Gr(A),\;(x',y')\in
Gr(B)\bigg\}.\end{equation}

The following lemmas are also needed. We refer to  \cite{KM} for
their proofs.
\begin{lem}\label{lem2.1}
Let $A$ be a maximal monotone operator of $\mathcal{H}$. If $x\in
\overline{D(A)}$ and  $y\in \mathcal{H}$ are such that
$$ \langle A^0(z)-y, z-x\rangle\geq 0\;\;\forall z\in D(A),$$
then $x\in D(A)$ and $y\in A(x)$.
\end{lem}
\begin{lem}\label{lem2.2}
Let $A_n$ $(n\in \mathbb{N})$ and $A$ be maximal monotone operators
of $\mathcal{H}$ such that $dis(A_n, A)\to 0$. Suppose also that
$x_n\in D(A_n)$ with $x_n\to x$ and $y_n\in A_n(x_n)$ with $y_n\to
y$ weakly for some $x, y\in E$. Then $x\in D(A)$ and $y\in A(x)$.
\end{lem}
\begin{lem}\label{lem2.3}
Let $A$ and $B$ be maximal monotone operators of $\mathcal{H}$. Then\\
1) for $\eta>0$ and $x\in D(A)$
$$ \|x-J_{\eta}^{B}(x)\| \leq \eta\|A^0(x)\|+dis(A,B)+\sqrt{\eta\big(1+\|A^0(x)\|\big)dis(A, B)}.$$
2) For $\eta>0$ and $x, x'\in E$
$$ \|J_{\eta}^{A}(x)-J_{\eta}^{A}(x')\|\leq
\|x-x'\|.$$
\end{lem}
\begin{lem}\label{lem2.4}
Let $A_n$ $(n\in \mathbb{N})$ and $A$ be maximal monotone operators
of $\mathcal{H}$ such that $dis(A_n, A)\to 0$ and $\|A^0_n(x)\|\leq
c(1+\|x\|)$ for some $c>0$, all $n\in \mathbb{N}$ and $x\in D(A_n)$.
Then for every $z\in D(A)$ there exists a sequence $(\zeta_n)$ such
that
 \begin{equation*}\label{2.2}\zeta_n\in
D(A_n),\;\;\;\zeta_n\to z\;\;\textmd{and}\;\;A_n^0(\zeta_n)\to
A^0(z).\end{equation*}
\end{lem}
\vskip2mm

We will establish our main result under the following hypotheses.\\
$(H_1)$ There exists a function $\rho:I\longrightarrow [0, +\infty[$
which is right continuous on $[0, T[$ and nondecreasing with
$\rho(0) = 0$ and $\rho(T)<+\infty$ such that
\begin{equation*}dis(A(t), A(s))\leq d\rho(]s,
t])=\rho(t)-\rho(s)\;\;\textmd{for}\;\;0\leq s\leq t\leq
T.\end{equation*} $(H_2)$ There exists a nonnegative real constant
$c$ such that
\begin{equation*}\|A^0(t,x)\|\leq c(1+\| x\|)\;\;\textmd{for}\;\;t\in I,\;x\in
D(A(t)).\end{equation*} $(H_3)$ There exists a nonnegative real
constant $m$ such that
$$ \|f(t,x)\|\leq m(1+\|x\|)\;\;\;\forall (t,x)\in I\times
\mathcal{H}.$$ $(H_4)$ There exists $\alpha\in L^1(I, \mathbb{R};
\lambda)$ such that
$$ \|f(t,x)-f(t,y)\|\leq \alpha(t)\|x-y\|\;\;\;\forall(t,x,y)\in I\times
\mathcal{H}\times \mathcal{H}.$$

\section{Main result}
Now we are able to state our main theorem.
\begin{thm}\label{Theorem 3.2}    Let for
every $t\in I$, $A(t):D(A(t))\subset \mathcal{H}\rightrightarrows
\mathcal{H}$ be a maximal monotone operator satisfying $(H_1)$ and
$(H_2)$. Let $f:I\times \mathcal{H}\longrightarrow \mathcal{H}$ be
such that for any fixed $x\in \mathcal{H}$, $f(\cdot, x)$ is
$\big(\mathcal{B}(I),\mathcal{B}(\mathcal{H})\big)$-measurable and
$(H_3)$, $(H_4)$ are satisfied.
\\Then for any   $u_0\in D(A(0))$, there exists a unique
BVRC solution $ u :I \longrightarrow \mathcal{H}$   to the problem
$$ (P_f)
\begin{cases}
-\displaystyle\frac{du}{d\nu }(t)\in A(t) u(t)+ f(t,u(t))
\frac{d\lambda} {d\nu}(t)  \;\;\;\nu-a.e.\,t\in
   I;\\
    u(t)\in D(A(t))\;\;\;\forall t\in I;\\
    u(0)=u_0.
\end{cases}
$$
\end{thm}

\proof Following  Castaing et al \cite{CM}, we can choose a sequence
$(\varepsilon_n)\subset ]0, 1]$ such that $\varepsilon_n \downarrow
0$ and a partition $0=t_0^n<t_1^n<\cdots<t_{q_n}^n=T$ of $I$, for
which we have
\begin{equation}\label{3.1}
\nu(]t^n_i, t^n _{i+1}])=|t_{i+1}^n-t_i^n|+d\rho(]t_i^n,
t_{i+1}^n])\leq
\varepsilon_n\;\;\;\textmd{for}\;i=0,...,q_n-1.\end{equation} For
each $i\in \{0,...,q_n-1\}$, put $I_i^n=]t_i^n, t_{i+1}^n]$,
\begin{equation}\label{3.2} \delta_{i+1}^n=d\rho(]t_i^n,
t_{i+1}^n]), \;\;\;\; \eta^n_{i+1} = t ^n_{i+1}-t ^n_i, \;\;\;\;
\beta ^n_{i+1}  = \nu(]t^n_i, t^n _{i+1}]).\end{equation} For every
$n\in \mathbb{N}$, put $u_0^n=u_0\in D(A(0))$, and let us define by
induction, the sequence $(u_i^n)_{0\leq i\leq q_n-1}$ such that
\begin{equation}\label{3.3} u_{i+1}^n= J^n_{i+1}
\Big(u_i^n-\int_{t_i^n}^{t_{i+1}^n} f(s,
u_i^n)d\lambda(s)\Big)\end{equation}  where
$J^n_{i+1}:=J^{A(t_{i+1}^n)}_{\beta^n_{i+1}} =\big(Id_{\mathcal{H}}+
\beta^n_{i+1} A(t_{i+1}^n)\big)^{-1}$.
\\  Clearly,  by the
properties of the resolvent,  we have  $u^n_{i+1}  \in
D(A(t_{i+1}^n))$ and
\begin{equation}\label{3.5} -\frac{1} {\beta^n_{i+1}}\Big(u^n_{i+1}-
u_i^n+\int_{t_i^n}^{t_{i+1}^n} f(s, u_i^n)d\lambda(s)\Big) \in A(
t_{i+1}^n) u^n _{i+1}. \end{equation} For $t\in [t^n_i, t
^n_{i+1}[$,  $i = 0,...,q_n-1$, set
\begin{equation}\label{3.7}
u_n(t)=   u_i^n + \frac{ \nu(]t^n_i, t])  }{ \nu(]t^n_i, t
^n_{i+1}])} \Big(u_{i+1}^n-u_i^n+\int_{t_i^n}^{t_{i+1}^n} f(s,u_i^n)
d\lambda(s)\Big)  - \int_{t_i^n}^{t}
f(s,u_i^n)d\lambda(s),\end{equation} and $u_n(T)=u_{q_n}^n$. So that
$u_n(t_i^n)=u_i^n$ and $u_n$ is a BVRC mapping on $I$.

\vskip2mm{\bf Step 1.}
 We prove in this step  that the sequence $(u_n)$ is bounded in  norm and
 in variation.\\Using Lemma \ref{lem2.3}, we have from
\eqref{3.3}, $(H_1)$,  $(H_2)$ and $(H_3)$, for $i=0,1,\cdots,
q_n-1$,
\begin{eqnarray*}
\| u_{i+1}^n-u_i^n\| &\leq& \Big\|J_{i+1}^n\Big(u_i^n-
\int_{t_i^n}^{t_{i+1}^n}
f(s,u_i^n)d\lambda(s)\Big)-J_{i+1}^n(u_i^n)\Big\|+\big\|J_{i+1}^n(u_i^n)-u_i^n\big\|\\
&\leq& \int_{t_i^n}^{t_{i+1}^n} \|f(s,u_i^n)\|
d\lambda(s)+\beta_{i+1}^n\|A^0(t_i^n,
u_i^n)\|+dis\big(A(t_{i+1}^n),A(t_i^n)\big)\\&+&\sqrt{\beta_{i+1}^n\big(1+\|A^0(t_i^n,
u_i^n)\|\big)dis\big(A(t_{i+1}^n),A(t_i^n)\big)}\\
&\leq& m(1+\|u_i^n\|)\beta_{i+1}^n
+\big(1+c(1+\|u_i^n\|)\big)\beta_{i+1}^n+
\sqrt{\big(1+c(1+\|u_i^n\|)\big)(\beta_{i+1}^n)^2}\\
&\leq& m(1+\|u_i^n\|)\beta_{i+1}^n
+\big(1+c(1+\|u_i^n\|)\big)\beta_{i+1}^n
+\big(1+c(1+\|u_i^n\|)\big)\beta_{i+1}^n,
\end{eqnarray*}
that is, \begin{equation}\label{3.8}\| u_{i+1}^n-u_i^n\|\leq \big(
(2c+m)\|u_i^n\|+2(1+c)+m\big)\beta_{i+1}^n,\end{equation} which
entails
$$ \|u_{i+1}^n\|\leq
\big(1+(2c+m)\beta_{i+1}^n\big)\|u_i^n\|+\big(2(1+c)+m\big)\
\beta_{i+1}^n.$$ Applying Lemma \ref{lem2.5}, we get
\begin{eqnarray*}\|u_i^n\|&\leq&
\Big(\|u_0\|+\big(2(1+c)+m\big)\sum_{j=0}^{i-1}
\beta_{j+1}^n\Big)\exp\Big((2c+m)\sum_{j=0}^{i-1}
\beta_{j+1}^n\Big)\\
&\leq&\Big(\|u_0\|+\big(2(1+c)+m\big)\nu(]0,T])\Big)\exp\Big((2c+m)
\nu(]0,T])\Big)=:m_1, \end{eqnarray*} using this last estimate in
\eqref{3.8}, we obtain
$$\|u_{i+1}^n-u_i^n\|\leq\Big((2c+m) m_1+2(1+c)+m\Big)\beta_{i+1}^n=:m_2\beta_{i+1}^n.$$
Consequently,
\begin{equation}\label{3.9}\|u_i^n\|\leq M\;\;\forall\, 0\leq i \leq
q_n\;\;\textmd{and}
 \;\;\|u_{i+1}^n-u_i^n\|\leq M\nu(]t_i^n,
t_{i+1}^n])\;\;\forall\, 0\leq i<q_n,
\end{equation}
where $M=\max(m_1, m_2)$. Now, for $t\in [t_i^n, t_{i+1}^n[$, we
have from \eqref{3.1}, \eqref{3.7} and  \eqref{3.9},
\begin{eqnarray}\label{3.10}
\|u_n(t)-u_i^n\| &=&  \Big\| \frac{ \nu(]t^n_i, t])}{ \nu(]t^n_i, t
^n_{i+1}])} \Big(u_{i+1}^n-u_i^n+\int_{t_i^n}^{t_{i+1}^n} f(s,
u_i^n) d\lambda(s)\Big)  - \int_{t_i^n}^{t}
f(s, u_i^n)d\lambda(s)\Big\|\nonumber\\
&\leq& \|u_{i+1}^n-u_i^n\|+2m(1+\|u_i^n\|)\delta_{i+1}^n \leq
M\varepsilon_n+2m(1+M) \varepsilon_n
=:M_1\varepsilon_n.\end{eqnarray} This last relation gives us, for
all $n\in \mathbb{N}$,
\begin{equation*}\|u_n(t)\|\leq M+M_1\varepsilon_n
\leq M+M_1=:M_2\;\;\;\forall t\in I,\end{equation*}  that is
\begin{equation}\label{3.11}
\sup_n\|u_n\|=\sup_n\big(\sup_{t\in I}\|u_n(t)\|\big)\leq M_2.
\end{equation}
%Now, if we fix
% $s\in [t_i^n, t_{i+1}^n[$ and $t\in [t_j^n, t_{j+1}^n[$ with
%$j>i$, we obtain by \eqref{3.1}, \eqref{3.9} and \eqref{3.10},
%\begin{eqnarray*}
%\|u_n(t)-u_n(s)\|&\leq&
%\|u_n(t)-u_j^n\|+\|u_j^n-u_i^n\|+\|u_n(s)-u_i^n\|\\
%&\leq&2M_1\varepsilon_n+ M \sum_{k=0}^{j-i-1} \nu\big(]t_{i+k}^n,
%t_{i+k+1}^n]\big)=2M_1\varepsilon_n+M\nu(]t_i^n, t_j^n])\\&\leq&
%2M_1\varepsilon_n+ M\nu(]t_i^n, t]) \leq 2M_1\varepsilon_n+
%M\big(\nu(]t_i^n, s])+\nu(]s, t])\big)\\&\leq&2M_1\varepsilon_n+
%M\big(\nu(]t_i^n, t_{i+1}^n[)+\nu(]s, t])\big)\\&\leq&
%2M_1\varepsilon_n+ M\varepsilon_n+M \nu(]s, t]).
%\end{eqnarray*}
%So that,  we get for $n\in \mathbb{N}$ and $0\leq s\leq t\leq T$,
%\begin{equation}\label{3.12} \|u_n(t)-u_n(s)\|\leq
%M\nu(]s,t])+(M+2M_1) \varepsilon_n.\end{equation}

{\bf Step 2.} Convergence of the sequences $(u_n)$ and $(\frac{du_n}{d\nu})$.\\
Define the functions $\phi_n,\;\theta_n: I\longrightarrow I$ by
$$\phi_n(t)=t_{i}^n,\;\;\;\theta_n(t)=t_{i+1}^n\;\;\textmd{for}\;\; t\in ]t_i^n,
t_{i+1}^n],\;i=0,1,...,q_n-1,\;\; \textmd{and}\;\;
\phi_n(0)=\theta_n(0)=0,
$$ and set for all $t\in I$,
\begin{equation*}
B_n(t)=\sum_{i=0}^{q_n-1}\frac{1}{\nu(]t_i^n, t_{i+1}^n])}\Big(
u_{i+1}^n-u_i^n+\int_{t_i^n}^{t_{i+1}^n}
f(s,u_i^n)d\lambda(s)\Big)\chi_{]t_i^n, t_{i+1}^n]}(t),
\end{equation*}
where $\chi_{J}(.)$ is the characteristic function of the set
$J\subset I$. Whence, by \eqref{3.7}, it is clear that
\begin{equation*}
u_n(t)=u_0+\int_{]0,t]} B_n(s)d\nu(s)-\int_{]0,t]} f(s,
u_n(\phi_n(s)))d\lambda(s).
\end{equation*}
Since $\frac{d\lambda}{d\nu}$ is a density of $\lambda$ w.r.t $\nu$,
 then by \eqref{7}, we get for every $t\in I$,
\begin{equation*}
u_n(t)=u_0+\int_{]0,t]} \Big(B_n(s)-f(s,
u_n(\phi_n(s)))\frac{d\lambda}{d\nu}(s)\Big)d\nu(s),
\end{equation*}
so that, $\frac{du_n}{d\nu}$ is a density of the vector measure
$du_n$ w.r.t $\nu$ and for $\nu$-almost every $t\in I$,
\begin{equation*}\label{8}\frac{du_n}{d\nu}(t)=B_n(t)-f(t,
u_n(\phi_n(t)))\frac{d\lambda}{d\nu}(t),\end{equation*} that is, by
the definition of $B_n$, for $\nu$-almost every $t\in I$,
\begin{equation}\label{8}
\frac{du_n}{d\nu}(t)+f(t,
u_n(\phi_n(t)))\frac{d\lambda}{d\nu}(t)=\sum_{i=0}^{q_n-1}\frac{1}{\nu(]t_i^n,
t_{i+1}^n])}\Big( u_{i+1}^n-u_i^n+\int_{t_i^n}^{t_{i+1}^n}
f(s,u_i^n)d\lambda(s)\Big)\chi_{]t_i^n, t_{i+1}^n]}(t).
\end{equation}
 From $(H_3)$, \eqref{3.9} and \eqref{8}, it results that
 \begin{equation}\label{9}
 \Big\|\frac{du_n}{d\nu}(t)+f(t,
u_n(\phi_n(t)))\frac{d\lambda}{d\nu}(t)\Big\|\leq
M+m(1+M)=:M_3\;\;\;\nu-a.e.\,t\in I.
 \end{equation}
On the other hand, we know that $\nu=d\rho+\lambda$, whence
$\frac{d\lambda}{d\nu}(t)\leq 1$. So that, by $(H_3)$ and
\eqref{3.9}, it follows that
\begin{equation}\label{10}
 \Big\|f(t,
u_n(\phi_n(t)))\frac{d\lambda}{d\nu}(t)\Big\|\leq
m(1+M)=:M_4\;\;\;\nu-a.e.\,t\in I,
 \end{equation}
 consequently, there is a Borel subset $J\subset I$ with $\nu(J)=0$,
 such that
 \begin{equation}\label{11}
 \Big\|\frac{du_n}{d\nu}(t)\Big\|\leq M_3+M_4=:M_5\;\;\;\forall t\in I\setminus J.
 \end{equation}
 Now, observe that by \eqref{3.5} and \eqref{8}, for each $n\in \mathbb{N}$, there is a
 Borel subset $J_n\subset I$ with $\nu(J_n)=0$ such that
 \begin{equation}\label{12}
 \frac{du_n}{d\nu}(t)+f(t,
u_n(\phi_n(t)))\frac{d\lambda}{d\nu}(t)\in
-A(\theta_n(t))u_n(\theta_n(t))\;\;\;\forall t\in I\setminus J_n,
\end{equation}
moreover,
\begin{equation}\label{13}
u_n(\theta_n(t))\in D\big(A(\theta_n(t))\big)\;\;\;\forall t\in I.
\end{equation}

We will show in the following, that $(u_n)$ is a Cauchy sequence.\\
Let $m, n\in \mathbb{N}$. Set for all $t\in I$,
$f_n(t):=f(t,u_n(\phi_n(t)))\frac{d\lambda}{d\nu}(t)$. From $(H_1)$
and relations \eqref{12}, \eqref{9}, we have for $\nu$-a.e. $t\in
I$,
\begin{eqnarray*}
&&\Big\langle u_n(\theta_n(t))-u_m(\theta_m(t)),
\frac{du_n}{d\nu}(t)+f_n(t)-\frac{du_m}{d\nu}(t)-f_m(t)\Big\rangle\\
&\leq&
\Big(1+\big\|\frac{du_n}{d\nu}(t)+f_n(t)\big\|+\big\|\frac{du_m}{d\nu}(t)
+f_m(t)\big\|\Big)
dis (A(\theta_n(t)), A(\theta_m(t)))\\
&\leq& (1+2M_3)|\rho(\theta_n(t))-\rho(\theta_m(t))|\leq
(1+2M_3)\big(|\rho(\theta_n(t))-t|+|\rho(\theta_m(t))-t|\big)\\
&\leq& (1+2M_3)\big(d\rho(]t, \theta_n(t)])+d\rho(]t,
\theta_m(t)])\big).
\end{eqnarray*}
On the other hand, from \eqref{9} and \eqref{11},
\begin{eqnarray*}
&&\Big\langle u_n(t)-u_n(\theta_n(t))-u_m(t)+u_m(\theta_m(t)),
\frac{du_n}{d\nu}(t)+f_n(t)-\frac{du_m}{d\nu}(t)-f_m(t)\Big\rangle\\
&\leq&
2M_3\big(\|u_n(t)-u_n(\theta_n(t))\|+\|u_m(t)-u_m(\theta_m(t))\|\big)\\
&\leq& 2M_3M_5\big(\nu(]t, \theta_n(t)])+\nu(]t, \theta_m(t)])\big).
\end{eqnarray*}
Combining these last relations, we get for $\nu$-a.e. $t\in I$,
\begin{eqnarray}\label{14}
&&\Big\langle u_n(t)-u_m(t),
\frac{du_n}{d\nu}(t)-\frac{du_m}{d\nu}(t)\Big\rangle\nonumber\\
&=& \Big\langle u_n(t)-u_n(\theta_n(t))-u_m(t)+u_m(\theta_m(t)),
\frac{du_n}{d\nu}(t)+f_n(t)-\frac{du_m}{d\nu}(t)-f_m(t)\Big\rangle\nonumber\\
&+& \Big\langle u_n(\theta_n(t))-u_m(\theta_m(t)),
\frac{du_n}{d\nu}(t)+f_n(t)-\frac{du_m}{d\nu}(t)-f_m(t)\Big\rangle\nonumber\\
&+&\Big\langle u_n(t)-u_m(t),-f_n(t)+f_m(t)\Big\rangle\nonumber\\
&\leq&  2M_3M_5\big(\nu(]t, \theta_n(t)])+\nu(]t, \theta_m(t)])\big)
+(1+2M_3)\big(d\rho(]t, \theta_n(t)])+d\rho(]t, \theta_m(t)])\big)\nonumber\\
&+&\Big\langle u_n(t)-u_m(t),-f_n(t)+f_m(t)\Big\rangle\nonumber\\
&\leq& c_1\big(\sigma_n(t)+\sigma_m(t)\big)+\Big\langle
u_n(t)-u_m(t),-f_n(t)+f_m(t)\Big\rangle,
\end{eqnarray}
where, $c_1:=2M_3M_5+(1+2M_3)$ and $\sigma_n(t)=\nu(]t,
\theta_n(t)])$. Let us examine the last term in the preceding
estimate. For all $t\in I$, we have by $(H_4)$
\begin{eqnarray*}
&&\Big\langle u_n(t)-u_m(t),-f_n(t)+f_m(t)\Big\rangle\\&=&
\Big\langle
u_n(t)-u_m(t),-f(t,u_n(\phi_n(t)))\frac{d\lambda}{d\nu}(t)+f(t,u_m(\phi_m(t)))\frac{d\lambda}{d\nu}(t)\Big\rangle\\
&\leq&\|u_n(t)-u_m(t)\|\big(\alpha(t)\frac{d\lambda}{d\nu}(t)\|u_n(\phi_n(t))-u_m(\phi_m(t))\|\big)\\
&\leq&\alpha(t)\frac{d\lambda}{d\nu}(t)\|u_n(t)-u_m(t)\|\big(\|u_n(\phi_n(t))-u_n(t)\|+\|u_n(t)-u_m(t)\|+\|u_m(t)-u_m(\phi_m(t))\|\big)\\
&\leq&
\alpha(t)\frac{d\lambda}{d\nu}(t)\|u_n(t)-u_m(t)\|^2+\gamma_{n,m}(t),
\end{eqnarray*}
where,
\begin{eqnarray*}&&\alpha(t)\frac{d\lambda}{d\nu}(t)\|u_n(t)-u_m(t)\|\big(\|u_n(\phi_n(t))-u_n(t)\|+\|u_m(t)-u_m(\phi_m(t))\|\big)\\
&\leq&M_5\alpha(t)\|u_n(t)-u_m(t)\|\big(\nu(]\phi_n(t),
t])+\nu(]\phi_m(t),
t])\big)\frac{d\lambda}{d\nu}(t)=:\gamma_{n,m}(t).
\end{eqnarray*}
Replacing in \eqref{14}, its results that for $\nu$-a.e. $t\in I$,
\begin{eqnarray}\label{15}
\Big\langle u_n(t)-u_m(t),
\frac{du_n}{d\nu}(t)-\frac{du_m}{d\nu}(t)\Big\rangle
&\leq&\alpha(t)\frac{d\lambda}{d\nu}(t)\|u_n(t)-u_m(t)\|^2+\Delta_{n,m}(t),
\end{eqnarray}
where $
\Delta_{n,m}(t):=c_1\big(\sigma_n(t)+\sigma_m(t)\big)+\gamma_{n,m}(t)$.
Since $\rho$ is right continuous, it is clear that for all $t\in I$,
$\sigma_n(t)\to 0$ as $n\to\infty$. On the other hand, we have for
all $t\in I$, $\nu(]\phi_n(t), t])\frac{d\lambda}{d\nu}(t)\to\nu(\{
t\})\frac{d\lambda}{d\nu}(t)=0$ thanks to \eqref{5}.  This shows
that $\gamma_{n,m}(t)\to 0$ as $n,m\to\infty$  since from
\eqref{3.11}, $\|u_n(t)\|\leq M_2$. We conclude that
\begin{equation}\label{16}\Delta_{n,m}(t)\longrightarrow 0\;\;\; \textmd{as}\;\; n,m\to\infty,\end{equation} and by the Lebesgue
dominated convergence theorem,
$$ \int_{]0, T]} \Delta_{n,m}(t)\,d\nu(t)\longrightarrow 0\;\;\;
\textmd{as}\;\;n,m\to\infty.$$ Applying Lemma \ref{Lem2.1}, we get
from \eqref{15}, taking account that $\|u_n(0)-u_m(0)\|=0$,
\begin{eqnarray*}
\|u_n(t)-u_m(t)\|^2\leq 2\int_{]0,t]}
\big(\alpha(s)\frac{d\lambda}{d\nu}(s)\|u_n(s)-u_m(s)\|^2\big)d\nu(s)+2\int_{]0,T]}\Delta_{n,m}(s)\,d\nu(s).
\end{eqnarray*}
Since (see \eqref{5})
$\alpha(s)\frac{d\lambda}{d\nu}(s)\nu(\{s\})=0<1$, we get by Lemma
\ref{lem2.6},
\begin{eqnarray*}
\|u_n(t)-u_m(t)\|^2&\leq&\Big(2\int_{]0,T]}\Delta_{n,m}(s)\,d\nu(s)\Big)\exp\Big(2\int_{]0,t]}
\big(\alpha(s)\frac{d\lambda}{d\nu}(s)\big)d\nu(s)\Big)\\&=&
\Big(2\int_{]0,T]}\Delta_{n,m}(s)\,d\nu(s)\Big)\exp\Big(2\int_{]0,t]}
\alpha(s)d\lambda(s)\Big),
\end{eqnarray*}
so that
$$\sup_{t\in I} \|u_n(t)-u_m(t)\|^2\leq\Big(2\int_{]0,T]}\Delta_{n,m}(s)\,d\nu(s)\Big)\exp\Big(2\int_{]0,T]}
\alpha(s)d\lambda(s)\Big).$$ We conclude, by \eqref{16}, that
$(u_n)$ is a Cauchy sequence w.r.t the norm of uniform convergence
on the space of $\mathcal{H}$-valued bounded mappings defined on
$I$. This implies that $(u_n)$ converges uniformly on $I$ to some
mapping $u$.

Now, observe that the sequence $(\frac{du_n}{d\nu})$ is bounded in
$L^2(I, \mathcal{H}; \nu)$ due to relation \eqref{11}, so that it
converges weakly in $L^2(I, \mathcal{H}; \nu)$ to some mapping $w\in
L^2(I, \mathcal{H}; \nu)$. In particular, for all $t\in I$
\begin{equation*}
\int_{]0, t]} \frac{du_n}{d\nu}(s)\, d\nu(s)\longrightarrow\int_{]0,
t]} w(s)\, d\nu(s)\;\;\;\textmd{weakly in }\;\mathcal{H}.
\end{equation*}
Since $\frac{du_n}{d\nu}$ is a density of the vector measure $du_n$
w.r.t the measure $\nu$, we have for all $t\in I$,
\begin{eqnarray*}
u_n(t)=u_0+\int_{]0, t]} \frac{du_n}{d\nu}(s)\,d\nu(s)
\end{eqnarray*}
and since $(u_n(t))$ converges strongly in $\mathcal{H}$ and then
weakly to $u(t)$, we deduce by what precedes that
\begin{equation*}
u(t)=u_0+\int_{]0, t]} w(s)\,d\nu(s)\;\;\;\forall t\in I.
\end{equation*}
This shows that $u$ is right continuous with bounded variation and
for $\nu$-a.e. $t\in I$, $\frac{du}{d\nu}=w$. Consequently,
$(\frac{du_n}{d\nu})$ converges weakly in $L^2(I, \mathcal{H}; \nu)$
to $\frac{du}{d\nu}$.

\vskip2mm

{\bf Step 3.} Existence of solution.\\
First, observe that from \eqref{13},  $ u_n(\theta_n(t))\in
D\big(A(\theta_n(t))\big)$ for all $t\in I $.  Also, we have that
for all $t\in I$,
$$dis\big(A(\theta_n(t)),A(t)\big)\leq |\rho(\theta_n(t))-\rho(t)|
\longrightarrow 0\;\;\;\textmd{as}\;\;n\to\infty,$$ and from
\eqref{11},
$$\|u_n(\theta_n(t))-u_n(t)\|\leq M_5\nu(]t, \theta_n(t)])\longrightarrow 0\;\;\;\textmd{as}\;\;n\to\infty,$$
so that, \begin{equation}\label{conv}\|u_n(\theta_n(t))-u(t)\|\leq
\|u_n(\theta_n(t))-u_n(t)\|+\|u_n(t)-u(t)\|\longrightarrow
0\;\;\;\textmd{as}\;\;n\to\infty.\end{equation} On the other hand,
by \eqref{3.11} and $(H_2)$, we have that the sequence
$\big(A^0\big(\theta_n(t), u_n(\theta_n(t))\big)\big)$ is bounded in
$\mathcal{H}$, and hence it is weakly relatively compact. Using all
these facts, we conclude by Lemma \ref{lem2.2}, that $u(t)\in
D(A(t))$ for all $t\in I$.

Next, since the sequence $(\frac{du_n}{d\nu})$ converges weakly in
$L^2(I, \mathcal{H}; \nu)$ to $\frac{du}{d\nu}$, by Mazur's theorem,
there is a sequence $(\xi_n)$ such that for each $n\in\mathbb{N}$,
$\xi_n\in co\{\frac{du_k}{d\nu};\;k\geq n\}$, and $(\xi_n)$
converges strongly in $L^2(I, \mathcal{H}; \nu)$ to
$\frac{du}{d\nu}$. So that, there is a subsequence $(\xi_{n_j})$,
which converges $\nu$-almost every where to $\frac{du}{d\nu}$. This
means the existence of a Borel subset $J'\subset I$, with
$\nu(J')=0$ and for $t\in I\setminus J'$,
$$\xi_{n_j}(t)\longrightarrow \frac{du}{d\nu}(t)\in\bigcap_{n}
\overline{co} \big\{\frac{du_{k}}{d\nu};\;k\geq n\big\}.$$ This
implies, by \eqref{co}, that for any fixed $\eta\in \mathcal{H}$,
\begin{equation}\label{17}
\big\langle \frac{du}{d\nu}(t), \eta
\big\rangle\leq\limsup_{n\to\infty}\big\langle \frac{du_n}{d\nu}(t),
\eta \big\rangle.
\end{equation}

To prove that $u$ is a solution to our considered problem, we will
use Lemma \ref{lem2.1}. Since $u(t)\in D(A(t))$ for all $t\in I$, we
have to show that for $\nu$-almost every fixed $t\in I$ and for any
$z\in D(A(t))$
\begin{equation*}
\big\langle
A^0(t,z)+\frac{du}{d\nu}(t)+f(t,u(t))\frac{d\lambda}{d\nu}(t),z-u(t)\big\rangle\geq
0.
\end{equation*}
Indeed, let $t\in I$. By $(H_2)$, using Lemma \ref{lem2.4}, we can
ensure the existence of a sequence $(\zeta_n)_n$, such that
\begin{equation}\label{18}
\zeta_n\in D(A(\theta_n(t))),\;\;\;\zeta_n\longrightarrow
z\;\;\textmd{and}\;\; A^0\big(\theta_n(t),
\zeta_n\big)\longrightarrow A^0(t,z).
\end{equation}
Since $A(\theta_n(t))$ is monotone, using \eqref{12}, we have for
each $n\in \mathbb{N}$ and $t\in I\setminus J_n$,
\begin{equation}\label{19}
\big\langle u_n(\theta_n(t))-\zeta_n,
\frac{du_n}{d\nu}(t)+f(t,u_n(\phi_n(t)))\frac{d\lambda}{d\nu}(t)+A^0(\theta_n(t),
\zeta_n)\big\rangle\leq 0.
\end{equation}
On the other hand, we have from $(H_4)$ and \eqref{11}
\begin{eqnarray*}
&&\|f(t,u_n(\phi_n(t))\frac{d\lambda}{d\nu}(t)-f(t,u(t))\frac{d\lambda}{d\nu}(t)\|\leq\alpha(t)\|u_n(\phi_n(t))-u(t)\|\frac{d\lambda}{d\nu}(t)\nonumber\\
&\leq&\alpha(t)(\|u_n(\phi_n(t))-u_n(t)\|+\|u_n(t)-u(t)\|)\frac{d\lambda}{d\nu}(t)\nonumber\\
&\leq&\alpha(t)(M_5\nu(]\phi_n(t),
t])+\|u_n(t)-u(t)\|)\frac{d\lambda}{d\nu}(t),
\end{eqnarray*}
whence by the uniform convergence of $(u_n)$ to $u$ and relation
\eqref{5}, it results that for all $t\in I$,
\begin{equation}\label{lip}
f(t,u_n(\phi_n(t))\frac{d\lambda}{d\nu}(t)\longrightarrow
f(t,u(t))\frac{d\lambda}{d\nu}(t)\;\;\;\textmd{as}\;\;n\to\infty.
\end{equation}

   Now,  let $t\in I\setminus \big(J\cup J'\cup(\underset{n}{\bigcup} J_n)\big)$. From \eqref{17}, we have
  \begin{equation*}
  \big\langle \frac{du}{d\nu}(t), u(t)-z\big\rangle\leq\limsup_{n\to\infty}\big\langle \frac{du_n}{d\nu}(t),
 u(t)-z\big\rangle,
  \end{equation*}
and since from  \eqref{11} and \eqref{19}, we have for each $n\in
\mathbb{N}$,
\begin{eqnarray*}
&&\big\langle \frac{du_n}{d\nu}(t),
  u(t)-z\big\rangle\leq\big\langle
  \frac{du_n}{d\nu}(t),u(t)-u_n(\theta_n(t))\big\rangle+ \big\langle\frac{du_n}{d\nu}(t), u_n(\theta_n(t))-\zeta_n\big\rangle+\big\langle
  \frac{du_n}{d\nu}(t),\zeta_n-z\big\rangle\\&\leq&
  M_5\big(\|u(t)-u_n(\theta_n(t))\|+\|\zeta_n-z\|\big)+\big\langle
  \frac{du_n}{d\nu}(t),u_n(\theta_n(t))-\zeta_n\big\rangle
  \\&\leq&
  M_5\big(\|u(t)-u_n(\theta_n(t))\|+\|\zeta_n-z\|\big)+\big\langle f(t,u_n(\phi_n(t)))\frac{d\lambda}{d\nu}(t)+A^0(\theta_n(t),
\zeta_n),\zeta_n-u_n(\theta_n(t))\big\rangle
\end{eqnarray*}
we conclude, using \eqref{conv}, \eqref{18} and \eqref{lip},
that\begin{equation} \big\langle \frac{du}{d\nu}(t),
u(t)-z\big\rangle\leq\big\langle
f(t,u(t))\frac{d\lambda}{d\nu}(t)+A^0(t, z),z-u(t)\big\rangle,
\end{equation}
which, with the fact that $u(t)\in D(A(t))$ and $u(0)=u_0$, means
that $u$ is a BVRC solution to our problem.

{\bf Step 4.} Uniqueness of the solution.\\
Let $u$, $v$ be two BVRC solutions of $(P_f)$. By the monotonicity
of the operator $A(t)$, we get for $\nu$-a.e. $t\in I$,
\begin{eqnarray*}
\big\langle
-\frac{du}{d\nu}(t)-f(t,u(t))\frac{d\lambda}{d\nu}(t)+\frac{dv}{d\nu}(t)+f(t,v(t))\frac{d\lambda}{d\nu}(t),
u(t)-v(t)\big\rangle\geq 0,
\end{eqnarray*}
that is
\begin{eqnarray*}
\big\langle \frac{du}{d\nu}(t)-\frac{dv}{d\nu}(t),
u(t)-v(t)\big\rangle\leq \big\langle
f(t,v(t))\frac{d\lambda}{d\nu}(t)-f(t,u(t))\frac{d\lambda}{d\nu}(t),
u(t)-v(t)\big\rangle,
\end{eqnarray*}
and by $(H_4)$, we obtain
\begin{eqnarray*}
\big\langle \frac{du}{d\nu}(t)-\frac{dv}{d\nu}(t),
u(t)-v(t)\big\rangle\leq\alpha(t)\frac{d\lambda}{d\nu}(t)\|u(t)-v(t)\|^2.
\end{eqnarray*}
Then, an application of Lemma \ref{Lem2.1} taking account that
$u(0)=v(0)=u_0$, gives us for all $t\in ]0, T]$,
\begin{eqnarray*}
\|u(t)-v(t)\|^2\leq 2\int_{]0, t]}
\alpha(s)\frac{d\lambda}{d\nu}(s)\|u(s)-v(s)\|^2 d\nu(s),
\end{eqnarray*}
from Lemma \ref{lem2.6}, it results that $\|u(t)-v(t)\|^2=0$ for all
$t\in I$, which means the uniqueness of the solution.
$\hfill\square$

\end{document}